\documentclass[11pt]{article}
\usepackage{amsfonts}
\usepackage{mathrsfs}
\usepackage[dvips]{graphics}
\usepackage[dvips]{color}
\usepackage{amsthm}
\usepackage[]{amsmath}
\usepackage{CJK}
\usepackage{indentfirst}
\newtheorem{proposition}{Proposition}[section]
\newtheorem{lemma}[proposition]{Lemma}

\newtheorem{theorem}[proposition]{Theorem}
\newtheorem{corollary}[proposition]{Corollary}
\newtheorem{property}[proposition]{Property}
\newtheorem{remark}[proposition]{Remark}
\newtheorem{algorithm}[proposition]{Algorithm}
\begin{document}
\begin{CJK*}{GBK}{song}
\CJKindent

\centerline{\textbf{\LARGE{The number of distinct and repeated squares}}}

\vspace{0.2cm}

\centerline{\textbf{\LARGE{and cubes in the Fibonacci sequence}}}

\vspace{0.2cm}

\centerline{Huang Yuke\footnote[1]{School of Mathematics and Systems Science, Beihang University (BUAA), Beijing, 100191, P. R. China. E-mail address: huangyuke07@tsinghua.org.cn,~hyg03ster@163.com.}
~~Wen Zhiying\footnote[2]{Department of Mathematical Sciences, Tsinghua University, Beijing, 100084, P. R. China. E-mail address: wenzy@tsinghua.edu.cn(Corresponding author).}}

\vspace{1cm}

\centerline{\textbf{\large{ABSTRACT}}}

\vspace{0.2cm}

The Fibonacci sequence $\mathbb{F}$ is the fixed point beginning with $a$ of morphism $\sigma(a,b)=(ab,a)$. In this paper, we get the explicit expressions of all squares and cubes, then we determine the number of distinct squares and cubes in $\mathbb{F}[1,n]$ for all $n$, where $\mathbb{F}[1,n]$ is the prefix of $\mathbb{F}$ of length $n$. By establishing and discussing the recursive structure
of squares and cubes, we give algorithms for counting the number of repeated squares and cubes in $\mathbb{F}[1,n]$ for all $n$, and get explicit expressions for some special $n$ such as $n=f_m$ (the Fibonacci number) etc., which including some known results such as in A.S.Fraenkel and J.Simpson\cite{FS1999,FS2014}, J.Shallit et al\cite{DMSS2014}.

\vspace{0.2cm}


\section{Introduction}

Let $\mathcal{A}=\{a,b\}$ be a binary alphabet.
The concatenation of factors $\nu$ and $\omega$ denoted by $\nu\omega$.
The Fibonacci sequence $\mathbb{F}$ is the fixed point beginning with $a$ of the Fibonacci morphism $\sigma$ defined by $\sigma(a)=ab$ and $\sigma(b)=a$.
As a classical example over a binary alphabet, $\mathbb{F}$ having many remarkable properties, we refer to M.Lothaire\cite{L1983,L2002}, J.M.Allouche and J.Shallit\cite{AS2003},
Berstel\cite{B1966,B1980}.

Let $\omega$ be a factor of $\mathbb{F}$, denoted by $\omega\prec\mathbb{F}$. Since $\mathbb{F}$ is uniformly recurrent, $\omega$ occurs infinitely many times.
Let $\omega_p$ be the $p$-th occurrence of $\omega$.
If the factor $\omega$ and integer $p$ such that $\omega_p\omega_{p+1}$ (resp. $\omega_p\omega_{p+1}\omega_{p+2}$) is the factor of $\mathbb{F}$, we call it a square (resp. cube) of $\mathbb{F}$.
As we know, $\mathbb{F}$ contains no fourth powers. The properties of squares and cubes are objects of a great interest in many aspects of mathematics and computer science etc.

We denote $F_m=\sigma^m(a)$ for $m\geq0$ and define $F_{-1}=b$, $F_{-2}=\varepsilon$ (empty word).
Define $f_m=|F_m|$ the $m$-th Fibonacci number,
$f_{-2}=0$, $f_{-1}=1$, $f_{m+1}=f_m+f_{m-1}$ for $m\geq -1$.
Let $\mathbb{F}[1,n]$ be the prefix of $\mathbb{F}$ of length $n$. In this paper, we consider the four functions below:

$A(n):=\sharp\{\omega:\omega\omega\prec\mathbb{F}[1,n]\}$, the number of distinct squares in $\mathbb{F}[1,n]$;

$B(n):=\sharp\{(\omega,p):\omega_p\omega_{p+1}\prec\mathbb{F}[1,n]\}$, the number of repeated squares in $\mathbb{F}[1,n]$;

$C(n):=\sharp\{\omega:\omega\omega\omega\prec\mathbb{F}[1,n]\}$, the number of distinct cubes in $\mathbb{F}[1,n]$;

$D(n):=\sharp\{(\omega,p):\omega_p\omega_{p+1}\omega_{p+2}\prec\mathbb{F}[1,n]\}$, the number of repeated cubes in $\mathbb{F}[1,n]$.

\vspace{0.2cm}

The methods for counting the four functions have attracted some many authors, but known results are not rich.
A.S.Fraenkel and J.Simpson gave the expression of $A(f_m)$ and $B(f_m)$ in 1999\cite{FS1999} and 2014\cite{FS2014}.
In 2014,
C.F.Du, H.Mousavi, L.Schaeffer and J.Shallit gave the expression of $B(f_m)$ and $D(f_m)$ by mechanical methods, see Theorem 58 and Theorem 59 in \cite{DMSS2014}.
In this paper,
we give the explicit expressions of $A(n)$, $B(f_m)$, $C(n)$ and $D(f_m)$.
Although we haven't get the explicit expressions of $B(n)$ and $D(n)$,
we give fast algorithms for counting $B(n)$ and $D(n)$ for all $n$.

\vspace{0.2cm}

The main tool of this paper is the ``structure properties" of the sequence of return words in the Fibonacci sequence, which introduced and studied in \cite{HW2015-1}, also see Property \ref{G}.
The definition of return words is from F.Durand\cite{D1998}. Let $\omega$ be a factor of $\mathbb{F}$. For $p\geq1$, let $\omega_p=x_{i+1}\cdots x_{i+n}$ and $\omega_{p+1}=x_{j+1}\cdots x_{j+n}$.
The factor $x_{i+1}\cdots x_{j}$ is called the $p$-th return word of $\omega$ and denoted by $r_{p}(\omega)$.
The sequence $\{r_p(\omega)\}_{p\geq1}$ is called the sequences of the return words of factor $\omega$.

By the ``structure properties" (Property \ref{G}), we can determine the positions of all $\omega_p$. By the definition of square (resp. cube) and return word, we have
$$\omega_p\omega_{p+1}\prec\mathbb{F}\Leftrightarrow r_p(\omega)=\omega,~~
\omega_p\omega_{p+1}\omega_{p+2}\prec\mathbb{F}\Leftrightarrow r_p(\omega)=r_{p+1}(\omega)=\omega,$$
where the ``=" means ``have the same expressions".
By these relations, we can determine the positions of all squares and cubes, and then get
$A(n)$, $B(n)$, $C(n)$ and $D(n)$. But this method is complicated, another improved and fast method is used in this paper.

\vspace{0.2cm}

This paper is organized as follows.
Section 2 present some basic notations and known results.
Section 3 prove some basic properties of squares.
We determine $A(n)$ (distinct squares) in Section 4.
Section 5 is devoted to establish the recursive structure of squares, then we
determine $B(n)$ (repeated squares) in
Section 6.
Similarly, we establish the recursive structure of cubes, then determine
$C(n)$ (distinct cubes) and $D(n)$ (repeated cubes) in Section 7 to 10.

\section{Preliminaries}

Let $\tau=x_1\cdots x_n$ be a finite word (or $\tau=x_1x_2\cdots$ be a sequence).
For any $i\leq j\leq n$, define $\tau[i,j]:=x_ix_{i+1}\cdots x_{j-1}x_j$.
By convention, we denote $\tau[i]:=\tau[i,i]=x_i$ and $\tau[i,i-1]=\varepsilon$.
Notation $\nu\triangleright\omega$ means word $\nu$ is a suffix of word $\omega$.

For $m\geq-1$, let $\delta_m\in\{a,b\}$ be the last letter of $F_m$, then
$\delta_m=a$ iff $m$ is even.
The $m$-th singular word is defined as $K_m=\delta_{m+1}F_m\delta_m^{-1}=\delta_{m+1}F_m[1,f_{m}-1]$ for $m\geq-1$.
By Property 2(9) in \cite{WW1994}, all singular words are palindromes.
Let $Ker(\omega)$ be the maximal singular word occurring in factor $\omega$,
called the kernel of $\omega$.
Then by Theorem 1.9 in \cite{HW2015-1},
$Ker(\omega)$ occurs in $\omega$ only once.
Moreover

\begin{property}[Theorem 2.8 in \cite{HW2015-1}]\label{wp}~
$Ker(\omega_p)=Ker(\omega)_p$ for all $\omega\in\mathbb{F}$ and $p\geq1.$
\end{property}
This means, let $Ker(\omega)=K_m$, then the maximal singular word occurring in $\omega_p$ is just $K_{m,p}$. For instance, $Ker(aba)=b$, $(aba)_3=\mathbb{F}[6,8]$, $(b)_3=\mathbb{F}[7]$, so $Ker((aba)_3)=(b)_3$, $(aba)_3=a(b)_3a$.

\begin{property}[Theorem 2.11 in \cite{HW2015-1}]\label{G}~
For any factor $\omega$,
the sequence of return words $\{r_p(\omega)\}_{p\geq1}$ is the Fibonacci sequence over the alphabet $\{r_1(\omega),r_2(\omega)\}$.
\end{property}

Property \ref{k1} and \ref{k2} are useful in our proofs. Property \ref{k1} can be proved by induction. Since all singular words are palindromes, Property \ref{k2} holds by the cylinder structure of palindromes in \cite{HW2016-1}.

\begin{property}[Lemma 2.2 in \cite{HW2015-1}]\label{k1}~For $m\geq-1$,
(1) $K_{m+3}=K_{m+1}K_{m}K_{m+1}$.

(2) $K_{m+2}=K_{m}K_{m+1}\delta_{m}^{-1}\delta_{m+1}=\delta_{m}^{-1}\delta_{m+1}K_{m+1}K_{m}$.
\end{property}

\begin{property}[]\label{k2}~$K_{m}\prec K_{m+3}[2,f_{m+3}-1]$, $K_{m+1}\not\!\prec K_{m+3}[2,f_{m+3}-1]$,
$K_{m+2}\not\!\prec K_{m+3}[2,f_{m+3}-1]$.
\end{property}

\section{Basic properties of squares}

By Definition 2.9 and Corollary 2.10 in \cite{HW2015-1}, any factor $\omega$ with kernel $K_m$ can be expressed uniquely as
$\omega=K_{m+1}[i,f_{m+1}] K_m K_{m+1}[1,j]=K_{m+3}[i,f_{m+2}+j],$
where $2\leq i\leq f_{m+1}+1$ and $0\leq j\leq f_{m+1}-1$.
By Property \ref{wp}, $\omega_p\omega_{p+1}\prec\mathbb{F}$ means
$$\omega_p\omega_{p+1}=K_{m+1}[i,f_{m+1}]
\underbrace{K_{m,p} K_{m+1}[1,j]K_{m+1}[i,f_{m+1}]}_{r_p(K_m)}
K_{m,p+1} K_{m+1}[1,j]\prec\mathbb{F}.$$

By Property \ref{G}, $K_m$ has only two distinct return words $r_1(K_m)=K_mK_{m+1}$ and $r_2(K_m)=K_mK_{m-1}$,
so $\omega_p\omega_{p+1}\prec\mathbb{F}$ has two cases as below,
and in each case, $|\omega|=|r_p(K_m)|$.

\vspace{0.2cm}

\textbf{Case 1.} $r_p(K_m)=r_1(K_m)=K_mK_{m+1}$. Comparing the two expressions of $r_p(K_m)$, we have
$$K_mK_{m+1}[1,j]K_{m+1}[i,f_{m+1}]=K_mK_{m+1}\Rightarrow
j=i-1.$$
Comparing the two ranges of $i$ that $2\leq i\leq f_{m+1}+1$ and $0\leq j=i-1\leq f_{m+1}-1$,
we have $2\leq i\leq f_{m+1}$ and $m\geq0$. Moreover $|\omega|=|r_1(K_m)|=f_{m+2}$ and
\begin{equation*}
\begin{split}
\omega\omega=&
K_{m+1}[i,f_{m+1}] K_{m} K_{m+1} K_{m} K_{m+1}[1,i-1]\\
=&K_{m+2}[i,f_{m+2}]\underline{K_{m+1}} K_{m+2}[1,f_{m}+i-1]
=K_{m+4}[i,2f_{m+2}+i-1].
\end{split}
\end{equation*}
The second and third equalities hold by Property \ref{k1}.

Since $K_{m+1}\prec\omega\omega\prec K_{m+4}[2,f_{m+4}-1]$, by Property \ref{k2}, $Ker(\omega\omega)=K_{m+1}$.

\vspace{0.2cm}

\textbf{Case 2.} $r_p(K_m)=r_2(K_m)=K_mK_{m-1}$.
Comparing the two expressions, we have
$$K_mK_{m+1}[1,j]K_{m+1}[i,f_{m+1}]=K_mK_{m-1}
\Rightarrow j=i-f_m-1.$$
So $f_m+1\leq i\leq f_{m+1}+1$ and $m\geq-1$.
Moreover $|\omega|=|r_2(K_m)|=f_{m+1}$ and
\begin{equation*}
\begin{split}
\omega\omega=&
K_{m+1}[i,f_{m+1}] K_{m} K_{m-1} K_{m} K_{m+1}[1,i-f_m-1]\\
=&K_{m+1}[i,f_{m+1}] \underline{K_{m+2}} K_{m+1}[1,i-f_m-1]
=K_{m+5}[f_{m+2}+i,f_{m+3}+f_{m+1}+i-1].
\end{split}
\end{equation*}

Since $K_{m+2}\prec\omega\omega\prec K_{m+5}[2,f_{m+4}-1]$, by Property \ref{k2}, $Ker(\omega\omega)=K_{m+2}$.

\begin{remark}
By the discussion above, we have:
all squares in $\mathbb{F}$ are of length $2f_m$ for some $m\geq0$; for all $m\geq0$,
there exists a square of length $2f_m$ in $\mathbb{F}$.
This is a known result of P.S$\acute{e}\acute{e}$bold\cite{S1985}.
\end{remark}

\begin{property}[Property 4.1 in \cite{HW2016-1}]\label{P}~
$P(K_m,p)=pf_{m+1}+(\lfloor\phi p\rfloor+1)f_{m}-1$ for $m\geq-1$, $p\geq1$.
\end{property}

\begin{corollary}[Corollary 4.2 in \cite{HW2016-1}]\label{P1}~
$P(a,p)=p+\lfloor\phi p\rfloor$,
$P(b,p)=2p+\lfloor\phi p\rfloor$ for $p\geq1$.
\end{corollary}

For $m,p\geq1$, we define two sets below
\begin{equation*}
\begin{cases}
\langle1,K_m,p\rangle:=
\{P(\omega\omega,p):Ker(\omega\omega)=K_m,|\omega|=f_{m+1},\omega\omega\prec\mathbb{F}\}\\
\langle2,K_m,p\rangle:=
\{P(\omega\omega,p):Ker(\omega\omega)=K_m,|\omega|=f_{m-1},\omega\omega\prec\mathbb{F}\}
\end{cases}
\end{equation*}
Obviously they correspond the two cases of squares respectively. By Property \ref{P} we have
\begin{equation*}
\begin{split}
\langle1,K_m,p\rangle=&\{P(\omega,p):
\omega=K_{m+1}[i,f_{m+1}]K_{m}K_{m+1}[1,f_{m-1}+i-1],2\leq i\leq f_{m}\}\\
=&\{P(K_m,p)+f_{m-1}+i-1,2\leq i\leq f_{m}\}\\
=&\{pf_{m+1}+\lfloor\phi p\rfloor f_{m}+f_{m+1},\cdots,pf_{m+1}+\lfloor\phi p\rfloor f_{m}+f_{m+2}-2\},\\
\langle2,K_m,p\rangle=&\{P(\omega,p):\omega=K_{m-1}[i,f_{m-1}] K_{m} K_{m-1}[1,i-f_{m-2}-1],f_{m-2}+1\leq i\leq f_{m-1}+1\}\\
=&\{P(K_m,p)+i-f_{m-2}-1,f_{m-2}+1\leq i\leq f_{m-1}+1\}\\
=&\{pf_{m+1}+\lfloor\phi p\rfloor f_{m}+f_{m}-1,\cdots,
pf_{m+1}+\lfloor\phi p\rfloor f_{m}+2f_{m-1}-1\}.
\end{split}
\end{equation*}

%

\begin{corollary} $\sharp\langle 1,K_m,p\rangle=f_m-1$ and $\sharp\langle 2,K_m,p\rangle=f_{m-3}+1$ for $m,p\geq1$.
\end{corollary}

\section{The number of distinct squares in $\mathbb{F}[1,n]$}

Denote $a(n):=\sharp\{\omega:\omega\omega\triangleright\mathbb{F}[1,n],
\omega\omega\not\!\prec\mathbb{F}[1,n-1]\}$,
obversely, $A(n)=\sum_{i=1}^n a(i)$.
In order to count $a(n)$, we only need to consider $\langle i,K_m,1\rangle$ where $i=1,2$.

\begin{property}[]
$\langle1,K_m,1\rangle=\{2f_{m+1},\cdots,f_{m+3}-2\}$,
$\langle2,K_m,1\rangle=\{f_{m+2}-1,\cdots,
f_{m+1}+2f_{m-1}-1\}$.
\end{property}

It is easy to see that sets $\langle i,K_m,1\rangle$ are pairwise disjoint, and each set contains some consecutive integers. Therefore we get a chain
$$\langle2,K_1,1\rangle,\langle1,K_1,1\rangle,
\cdots,
\langle1,K_{m-1},1\rangle,\langle2,K_m,1\rangle,
\langle1,K_m,1\rangle,\langle2,K_{m+1},1\rangle,\cdots$$

By this chain, $a(n)=1$ iff $n\in\cup_{m\geq1}(\langle2,K_m,1\rangle\cup\langle1,K_m,1\rangle)$.
The ``$\cup$" means pairwise disjoint union in this paper.
Moreover, we have $\langle1,K_m,1\rangle\cup\langle2,K_{m+1},1\rangle=\{2f_{m+1},\cdots,f_{m+2}+2f_{m}-1\}$.

\begin{property}[] $a(1)=a(2)=a(3)=0$, $a(4)=1$ and for $n\geq5$
$$a(n)=1\text{ iff }n\in\cup_{m\geq1}\{2f_{m+1},\cdots,f_{m+2}+2f_{m}-1\}.$$
\end{property}

One method for counting $A(n)$ is by $A(n)=\sum_{i=1}^n a(i)$. By consider $A(f_{m+2}+2f_{m}-1)$ for $m\geq1$, we can give a fast algorithm of $A(n)$ for all $n\geq1$.
Since $\sum_{i=-1}^mf_i= f_{m+2}-1$,
$$\begin{array}{rl}
&A(f_{m+2}+2f_{m}-1)=a(4)+\sum\limits_{i=1}^{m}\sharp\{2f_{i+1},\cdots,f_{i+2}+2f_{i}-1\}\\
=&1+\sum\limits_{i=1}^{m}(f_{i}+f_{i-2})
=1+\sum\limits_{i=-1}^{m}f_{i}-f_0-f_{-1}+\sum\limits_{i=-1}^{m-2}f_{i}
=f_{m+2}+f_{m}-3.
\end{array}$$

\begin{theorem}[]\label{T4.3} For all $n\geq1$, let $m$ satisfies $2f_m\leq n<2f_{m+1}$,
\begin{equation*}
A(n)=\begin{cases}
n-f_{m-1}-2,&n\leq f_{m+1}+2f_{m-1}-1;\\
f_{m+1}+f_{m-1}-3,&otherwise.
\end{cases}
\end{equation*}
\end{theorem}

\begin{proof} When $f_{m+1}+2f_{m-1}\leq n\leq 2f_{m+1}-1$, $a(n)=0$,
$A(n)=A(f_{m+1}+2f_{m-1}-1)=f_{m+1}+f_{m-1}-3$.

When $2f_m\leq n\leq f_{m+1}+2f_{m-1}-1$, $a(n)=1$, $A(n)=A(2f_m-1)+n-2f_m+1$. Since $A(2f_m-1)=A(f_{m}+2f_{m-2}-1)$, we have $A(n)=n-f_{m-1}-2$. Thus the conclusion holds.
\end{proof}

%
%

\begin{remark}
Since $2f_{m-2}\leq f_{m}\leq f_{m-1}+2f_{m-3}-1$ for $m\geq2$,
as a spacial case of Theorem \ref{T4.3},
$$A(f_{m})=f_{m}-f_{m-3}-2=2f_{m-2}-2.$$
This is a known result of  A.S.Fraenkel and J.Simpson, see Theorem 1 in \cite{FS1999}.
\end{remark}

\section{The recursive structure of squares}

In this section, we establish a recursive structure of squares. Using it, we will count
the number of repeated squares in $\mathbb{F}[1,n]$ (i.e. $B(n)$) in Section 6.
For $m,p\geq1$, consider the vectors
\begin{equation*}
\begin{split}
\Gamma_{1,m,p}:&=[pf_{m+1}+\lfloor\phi p\rfloor f_{m}+f_{m+1}-1,\langle 1,K_m,p\rangle]\\
&=[pf_{m+1}+\lfloor\phi p\rfloor f_{m}+f_{m+1}-1,\cdots,pf_{m+1}+\lfloor\phi p\rfloor f_{m}+f_{m+2}-2];\\
\Gamma_{2,m,p}:&=[\langle 2,K_m,p\rangle,pf_{m+1}+\lfloor\phi p\rfloor f_{m}+2f_{m-1},\cdots,pf_{m+1}+\lfloor\phi p\rfloor f_{m}+f_{m+1}-2]\\
&=[pf_{m+1}+\lfloor\phi p\rfloor f_{m}+f_{m}-1,\cdots,pf_{m+1}+\lfloor\phi p\rfloor f_{m}+f_{m+1}-2].
\end{split}
\end{equation*}
Here vector $[\langle i,K_m,p\rangle]$ means arrange all elements in set $\langle i,K_m,p\rangle$, $i=1,2$.

Obversely, each $\Gamma_{i,m,p}$ contains consecutive integers.
The numbers of components in vectors $\Gamma_{1,m,p}$ and $\Gamma_{2,m,p}$ are
$f_{m}$ and $f_{m-1}$ respectively.
Moreover $\max\Gamma_{2,m,p}+1=\min\Gamma_{1,m,p}$ for $m,p\geq1$.


\begin{lemma}[Lemma 5.3 and 5.4 in \cite{HW2016-1}]\label{L}~
$\lfloor\phi(p+\lfloor\phi p\rfloor+1)\rfloor=p$,
$\lfloor\phi(2p+\lfloor\phi p\rfloor+1)\rfloor=p+\lfloor\phi p\rfloor$.
\end{lemma}

\begin{property}[]\label{R1} 
$\Gamma_{1,m,p}=[\Gamma_{2,m-1,P(a,p)+1},\Gamma_{1,m-1,P(a,p)+1}]$ for $m\geq2$, $p\geq1$.
\end{property}

\begin{proof}
By Corollary \ref{P1}, $P(a,p)+1=p+\lfloor\phi p\rfloor+1$. By Lemma \ref{L}, $\lfloor\phi (p+\lfloor\phi p\rfloor+1)\rfloor=p$.
\begin{equation*}
\begin{split}
&\min\Gamma_{2,m-1,P(a,p)+1}
=(p+\lfloor\phi p\rfloor+1)f_{m}+\lfloor\phi(p+\lfloor\phi p\rfloor+1)\rfloor f_{m-1}+f_{m-1}-1\\
=&(p+\lfloor\phi p\rfloor+1)f_{m}+p f_{m-1}+f_{m-1}-1=pf_{m+1}+\lfloor\phi p\rfloor f_{m}+f_{m+1}-1=\min\Gamma_{1,m,p};\\
&\max\Gamma_{1,m-1,P(a,p)+1}
=(p+\lfloor\phi p\rfloor+1)f_{m}+\lfloor\phi(p+\lfloor\phi p\rfloor+1)\rfloor f_{m-1}+f_{m+1}-2\\
=&(p+\lfloor\phi p\rfloor+1)f_{m}+p f_{m-1}+f_{m+1}-2=pf_{m+1}+\lfloor\phi p\rfloor f_{m}+f_{m+2}-2=\max\Gamma_{1,m,p}.
\end{split}
\end{equation*}

Since $\max\Gamma_{2,m,p}+1=\min\Gamma_{1,m,p}$ for $m,p\geq1$,
$\max\Gamma_{2,m-1,P(a,p)+1}+1=\min\Gamma_{1,m-1,P(a,p)+1}$.
Thus the conclusion holds.
\end{proof}

By an analogous argument, we have

\begin{property}[]\label{R2}
$\Gamma_{2,m,p}=[\Gamma_{2,m-2,P(b,p)+1},\Gamma_{1,m-2,P(b,p)+1}]$ for $m\geq3$, $p\geq1$.
\end{property}

In Property \ref{R1} and \ref{R2}, we establish the recursive relations for any
$\Gamma_{1,m,p}$ ($m\geq2$) and $\Gamma_{2,m,p}$ ($m\geq3$).
By the one-to-one correspondence between $\Gamma_{i,m,p}$ and $\langle i,K_m,p\rangle$, we can define the recursive structure over $\{\langle i,K_m,p\rangle|~i=1,2;~m,p\geq1\}$ denoted by $\mathcal{S}$. Each $\langle i,K_m,p\rangle$ is an element in $\mathcal{S}$.
The recursive structure $\mathcal{S}$ is a family of finite trees with root $\langle i,K_m,1\rangle$ for all $i=1,2$, $m\geq1$; and with recursive relations:
\begin{equation*}
\begin{cases}
\tau_1\langle1,K_m,p\rangle=\langle2,K_{m-1},P(a,p)+1\rangle\cup\langle1,K_{m-1},P(a,p)+1\rangle&\text{for }m\geq2;\\
\tau_2\langle2,K_m,p\rangle=\langle2,K_{m-2},P(b,p)+1\rangle\cup\langle1,K_{m-2},P(b,p)+1\rangle&\text{for }m\geq3.
\end{cases}
\end{equation*}

\begin{property}[]\
Each $\langle i,K_m,p\rangle$ belongs to the recursive structure $\mathcal{S}$,
$i=1,2$, $m,p\geq1$.
\end{property}

\begin{proof} Each element $\langle i,K_m,1\rangle$ is root of a finite tree in $\mathcal{S}$.
For $m,p\geq1$,
\begin{equation*}
\begin{cases}
\langle1,K_m,P(a,p)+1\rangle\in\tau_1\langle1,K_{m+1},p\rangle\\
\langle1,K_m,P(b,p)+1\rangle\in\tau_2\langle2,K_{m+2},p\rangle
\end{cases}
\text{and }
\begin{cases}
\langle2,K_m,P(a,p)+1\rangle\in\tau_1\langle1,K_{m+1},p\rangle\\
\langle2,K_m,P(b,p)+1\rangle\in\tau_2\langle2,K_{m+2},p\rangle
\end{cases}
\end{equation*}
Since $\mathbb{N}=\{1\}\cup\{P(a,p)+1\}\cup\{P(b,p)+1\}$, the recursive structure $\mathcal{S}$ contains all $\langle i,K_m,p\rangle$.
\end{proof}

On the other hand, by the recursive relations $\tau_1$ and $\tau_2$, each element $\langle i,K_m,p\rangle$ has a unique position in $\mathcal{S}$.
By Property \ref{R1} and \ref{R2}, the trees in $\mathcal{S}$ are pairwise disjoint.
Fig.1 and Fig.2 show the two finite trees in the recursive structure $\mathcal{S}$ with roots $\langle1,K_5,1\rangle$ and $\langle 2,K_5,1\rangle$ respectively.

\scriptsize
\setlength{\unitlength}{0.85mm}
\begin{center}
\begin{picture}(155,56)
\put(126,1){53}
\put(120.5,5){$\langle1,K_1,12\rangle$}
\put(120,0){\line(1,0){16}}
\put(120,8){\line(1,0){16}}
\put(120,0){\line(0,1){8}}
\put(136,0){\line(0,1){8}}
\put(146,9){51}
\put(140.5,13){$\langle2,K_1,12\rangle$}
\put(140,8){\line(1,0){16}}
\put(140,16){\line(1,0){16}}
\put(140,8){\line(0,1){8}}
\put(156,8){\line(0,1){8}}
\put(95,1){53}
\put(95,5){52}
\put(90.5,9){$\langle1,K_2,7\rangle$}
\put(90,0){\line(1,0){15}}
\put(90,12){\line(1,0){15}}
\put(90,0){\line(0,1){12}}
\put(105,0){\line(0,1){12}}
\put(110,13){50}
\put(110,17){49}
\put(105.5,21){$\langle2,K_2,7\rangle$}
\put(105.5,12.5){\line(1,0){14}}
\put(105.5,24){\line(1,0){14}}
\put(105.5,12.5){\line(0,1){11.5}}
\put(119.5,12.5){\line(0,1){11.5}}
\put(65,1){53}
\put(65,5){52}
\put(65,9){51}
\put(65,13){50}
\put(60.5,17){$\langle1,K_3,4\rangle$}
\put(60,0){\line(1,0){15}}
\put(60,20){\line(1,0){15}}
\put(60,0){\line(0,1){20}}
\put(75,0){\line(0,1){20}}
\put(35,1){53}
\put(35,5){52}
\put(35,9){51}
\put(35,13){50}
\put(35,17){49}
\put(35,21){48}
\put(35,25){47}
\put(30.5,29){$\langle1,K_4,2\rangle$}
\put(30,0){\line(1,0){15}}
\put(30,32){\line(1,0){15}}
\put(30,0){\line(0,1){32}}
\put(45,0){\line(0,1){32}}
\put(126,21){48}
\put(120.5,25){$\langle1,K_1,11\rangle$}
\put(120.5,20){\line(1,0){15.5}}
\put(120.5,28){\line(1,0){15.5}}
\put(120.5,20){\line(0,1){8}}
\put(136,20){\line(0,1){8}}
\put(80,25){47}
\put(80,29){46}
\put(75.5,33){$\langle2,K_3,4\rangle$}
\put(75,24){\line(1,0){14.5}}
\put(75,36){\line(1,0){14.5}}
\put(75,24){\line(0,1){12}}
\put(89.5,24){\line(0,1){12}}
\put(146,29){46}
\put(140.5,33){$\langle2,K_1,11\rangle$}
\put(140,28){\line(1,0){16}}
\put(140,36){\line(1,0){16}}
\put(140,28){\line(0,1){8}}
\put(156,28){\line(0,1){8}}
\put(95,33){45}
\put(95,37){44}
\put(90.5,41){$\langle1,K_2,6\rangle$}
\put(90.5,32){\line(1,0){14.5}}
\put(90.5,44){\line(1,0){14.5}}
\put(90.5,32){\line(0,1){12}}
\put(105,32){\line(0,1){12}}
\put(146,41){43}
\put(140.5,45){$\langle2,K_1,10\rangle$}
\put(140,40){\line(1,0){16}}
\put(140,48){\line(1,0){16}}
\put(140,40){\line(0,1){8}}
\put(156,40){\line(0,1){8}}
\put(126,33){45}
\put(120.5,37){$\langle1,K_1,10\rangle$}
\put(120,32){\line(1,0){16}}
\put(120,40){\line(1,0){16}}
\put(120,32){\line(0,1){8}}
\put(136,32){\line(0,1){8}}
\put(50,41){43}
\put(50,45){42}
\put(50,49){41}
\put(45.5,53){$\langle2,K_4,2\rangle$}
\put(45,40){\line(1,0){15}}
\put(45,56){\line(1,0){15}}
\put(45,40){\line(0,1){16}}
\put(60,40){\line(0,1){16}}
\put(110,45){42}
\put(110,49){41}
\put(105.5,53){$\langle2,K_2,6\rangle$}
\put(105.5,44){\line(1,0){14.5}}
\put(105.5,56){\line(1,0){14.5}}
\put(105.5,44){\line(0,1){12}}
\put(120,44){\line(0,1){12}}
\put(5,1){53}
\put(5,5){52}
\put(5,9){51}
\put(5,13){50}
\put(5,17){49}
\put(5,21){48}
\put(5,25){47}
\put(5,29){46}
\put(5,33){45}
\put(5,37){44}
\put(5,41){43}
\put(5,45){42}
\put(0.5,49){$\langle1,K_5,1\rangle$}
\put(0,0){\line(1,0){15}}
\put(0,52){\line(1,0){15}}
\put(0,0){\line(0,1){52}}
\put(15,0){\line(0,1){52}}
\put(63,52){$\tau_2$}
\put(61,50){\vector(3,-1){28}}
\put(61,50){\vector(1,0){43}}
\put(18,38){$\tau_1$}
\put(16,34){\vector(1,-1){13}}
\put(16,34){\vector(2,1){28}}
\put(108,37){$\tau_1$}
\put(106,42){\vector(2,-1){13}}
\put(106,42){\vector(1,0){33}}
\put(93,29.5){$\tau_2$}
\put(91,26){\vector(1,0){28}}
\put(91,26){\line(3,1){15}}
\put(106,31){\vector(1,0){33}}
\put(48,29){$\tau_1$}
\put(46,27){\vector(1,-1){13}}
\put(46,27){\vector(1,0){28}}
\put(78,19){$\tau_1$}
\put(76,17){\vector(1,-1){13}}
\put(76,17){\vector(1,0){28}}
\put(108,4){$\tau_1$}
\put(106,10){\vector(2,-1){13}}
\put(106,10){\vector(1,0){33}}
\end{picture}
\end{center}
\normalsize
\vspace{-0.2cm}
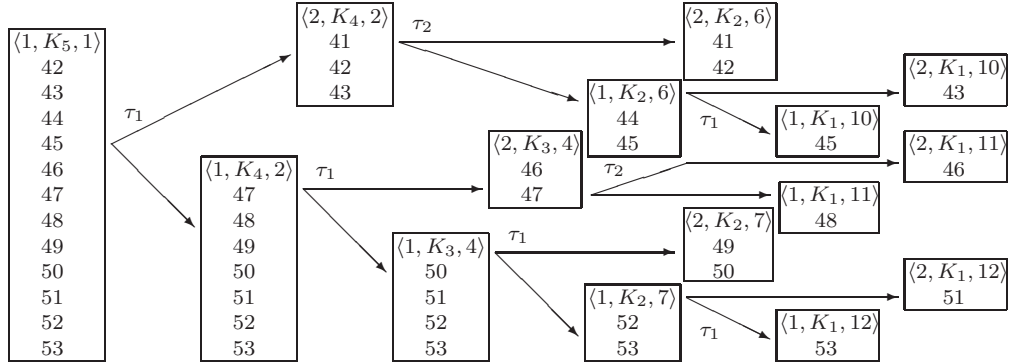
\centerline{Fig.1: The finite tree in the recursive structure $\mathcal{S}$ with root $\langle1,K_5,1\rangle$.}

\scriptsize
\setlength{\unitlength}{0.85mm}
\begin{center}
\begin{picture}(150,40)
\put(5,17){36}
\put(5,21){35}
\put(5,25){34}
\put(5,29){33}
\put(0.5,33){$\langle2,K_5,1\rangle$}
\put(0,16){\line(1,0){15}}
\put(0,36){\line(1,0){15}}
\put(0,16){\line(0,1){20}}
\put(15,16){\line(0,1){20}}
\put(19,24){$\tau_2$}
\put(16,30){\vector(2,-1){28}}
\put(16,30){\vector(1,0){43}}
\put(50,1){40}
\put(50,5){39}
\put(50,9){38}
\put(50,13){37}
\put(45.5,17){$\langle1,K_3,3\rangle$}
\put(45,20){\line(1,0){15}}
\put(45,0){\line(1,0){15}}
\put(45,0){\line(0,1){20}}
\put(60,0){\line(0,1){20}}
\put(64,13){$\tau_1$}
\put(61,10){\vector(2,-1){13}}
\put(61,10){\vector(3,1){28}}
\put(65,25){34}
\put(65,29){33}
\put(60.5,33){$\langle2,K_3,3\rangle$}
\put(60,24){\line(1,0){15}}
\put(60,36){\line(1,0){15}}
\put(60,24){\line(0,1){12}}
\put(75,24){\line(0,1){12}}
\put(79,37){$\tau_2$}
\put(76,35){\vector(3,-1){30}}
\put(76,35){\vector(1,0){48}}
\put(80,1){40}
\put(80,5){39}
\put(75.5,9){$\langle1,K_2,5\rangle$}
\put(75,0){\line(1,0){15}}
\put(75,12){\line(1,0){15}}
\put(75,0){\line(0,1){12}}
\put(90,0){\line(0,1){12}}
\put(94,1){$\tau_1$}
\put(91,4){\vector(1,0){15}}
\put(91,4){\vector(3,1){33}}
\put(95,13){37}
\put(95,17){36}
\put(90.5,21){$\langle2,K_2,5\rangle$}
\put(90,12.5){\line(1,0){15}}
\put(90,24){\line(1,0){15}}
\put(90,12.5){\line(0,1){11.5}}
\put(105,12.5){\line(0,1){11.5}}
\put(112,1){40}
\put(107.5,5){$\langle1,K_1,9\rangle$}
\put(107,0){\line(1,0){15}}
\put(107,8){\line(1,0){15}}
\put(107,0){\line(0,1){8}}
\put(122,0){\line(0,1){8}}
\put(112,21){35}
\put(107.5,25){$\langle1,K_1,8\rangle$}
\put(107,20){\line(1,0){15}}
\put(107,28){\line(1,0){15}}
\put(107,20){\line(0,1){8}}
\put(122,20){\line(0,1){8}}
\put(130,9){38}
\put(125.5,13){$\langle2,K_1,9\rangle$}
\put(125,8){\line(1,0){15}}
\put(125,16){\line(1,0){15}}
\put(125,8){\line(0,1){8}}
\put(140,8){\line(0,1){8}}
\put(130,29){33}
\put(125.5,33){$\langle2,K_1,8\rangle$}
\put(125,28){\line(1,0){15}}
\put(125,36){\line(1,0){15}}
\put(125,28){\line(0,1){8}}
\put(140,28){\line(0,1){8}}
\end{picture}
\end{center}
\normalsize
\vspace{-0.2cm}
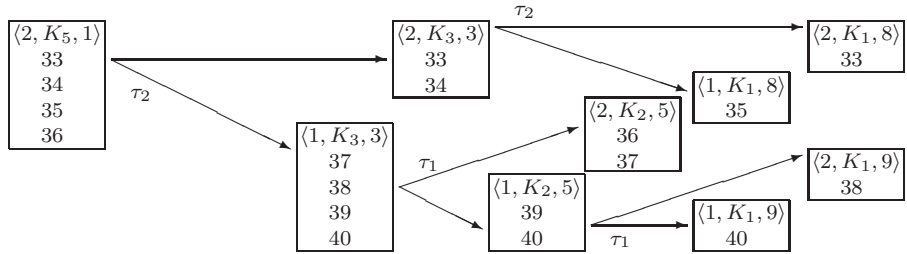
\centerline{Fig.2: The finite tree in the recursive structure $\mathcal{S}$ with root $\langle2,K_5,1\rangle$.}

\vspace{0.2cm}

By the recursive structure $\mathcal{S}$, we have the relation between the number of squares ending at position $\Gamma_{1,m,p}[i]=pf_{m+1}+\lfloor\phi p\rfloor f_{m}+f_{m+1}+i-1$ and $\Gamma_{1,m,1}[i]=2f_{m+1}+i-1$, see Property \ref{P5.5}.
Similarly, we have the relation between the number of squares ending at position $\Gamma_{2,m,p}[i]=pf_{m+1}+\lfloor\phi p\rfloor f_{m}+f_{m}+i-2$ and $\Gamma_{2,m,1}[i]=f_{m+2}+i-2$, see Property \ref{P5.6}.

\begin{property}[]\label{P5.5} For $1\leq i\leq f_m-1$,
\begin{equation*}
\begin{split}
&\{\omega:\omega\omega\triangleright\mathbb{F}[1,2f_{m+1}+i-1]\}\\
=&\{\omega:\omega\omega\triangleright\mathbb{F}[1,pf_{m+1}+\lfloor\phi p\rfloor f_{m}+f_{m+1}+i-1],Ker(\omega)=K_j,1\leq j\leq m\}.
\end{split}
\end{equation*}
\end{property}

\begin{property}[]\label{P5.6} For $1\leq i\leq f_{m-3}+1$,
\begin{equation*}
\begin{split}
&\{\omega:\omega\omega\triangleright\mathbb{F}[1,f_{m+2}+i-2]\}\\
=&\{\omega:\omega\omega\triangleright\mathbb{F}[1,pf_{m+1}+\lfloor\phi p\rfloor f_{m}+f_{m}+i-2],Ker(\omega)=K_j,1\leq j\leq m\}.
\end{split}
\end{equation*}
\end{property}

For instance, taking $m=3$, $p=3$ and $i=2$ in the property above. All squares ending at position 13 are $\{aabaab\}$.
All squares ending at position
34 are $\{aabaab,abaabaababaabaab\}$.
Since $Ker(aabaab)=aabaa=K_{3}$ and $Ker(abaabaababaabaab)=aabaababaabaa=K_5$,
only $\{aabaab\}$ is square with kernel $K_{j}$, $1\leq j\leq 3$.
Fig.3 shows the relation:
\scriptsize
\setlength{\unitlength}{0.85mm}
\begin{center}
\begin{picture}(150,48)
\put(5,2){36}
\put(5,6){35}
\put(5,10){34}
\put(5,14){33}
\put(0.5,18){$\langle2,K_5,1\rangle$}
\put(0,1){\line(1,0){15}}
\put(0,21){\line(1,0){15}}
\put(0,1){\line(0,1){20}}
\put(15,1){\line(0,1){20}}
\put(19,9){$\tau_2$}
\put(16,15){\vector(2,-1){28}}
\put(16,15){\vector(1,0){43}}
\put(48,0){$\cdots\cdots$}
\put(65,10){34}
\put(65,14){33}
\put(60.5,18){$\langle2,K_3,3\rangle$}
\put(60,9){\line(1,0){15}}
\put(60,21){\line(1,0){15}}
\put(60,9){\line(0,1){12}}
\put(75,9){\line(0,1){12}}
\put(79,22){$\tau_2$}
\put(76,20){\vector(3,-1){30}}
\put(76,20){\vector(1,0){48}}
\put(112,6){35}
\put(107.5,10){$\langle1,K_1,8\rangle$}
\put(107,5){\line(1,0){15}}
\put(107,13){\line(1,0){15}}
\put(107,5){\line(0,1){8}}
\put(122,5){\line(0,1){8}}
\put(130,14){33}
\put(125.5,18){$\langle2,K_1,8\rangle$}
\put(125,13){\line(1,0){15}}
\put(125,21){\line(1,0){15}}
\put(125,13){\line(0,1){8}}
\put(140,13){\line(0,1){8}}
\put(67,24){$\vdots$}
\put(114,24){$\vdots$}
\put(132,24){$\vdots$}
\put(65,35){13}
\put(65,39){12}
\put(60.5,43){$\langle2,K_3,1\rangle$}
\put(60,34){\line(1,0){15}}
\put(60,46){\line(1,0){15}}
\put(60,34){\line(0,1){12}}
\put(75,34){\line(0,1){12}}
\put(79,47){$\tau_2$}
\put(76,45){\vector(3,-1){30}}
\put(76,45){\vector(1,0){48}}
\put(112,31){14}
\put(107.5,35){$\langle1,K_1,3\rangle$}
\put(107,30){\line(1,0){15}}
\put(107,38){\line(1,0){15}}
\put(107,30){\line(0,1){8}}
\put(122,30){\line(0,1){8}}
\put(130,39){12}
\put(125.5,43){$\langle2,K_1,3\rangle$}
\put(125,38){\line(1,0){15}}
\put(125,46){\line(1,0){15}}
\put(125,38){\line(0,1){8}}
\put(140,38){\line(0,1){8}}
\end{picture}
\end{center}
\normalsize
\vspace{-0.2cm}
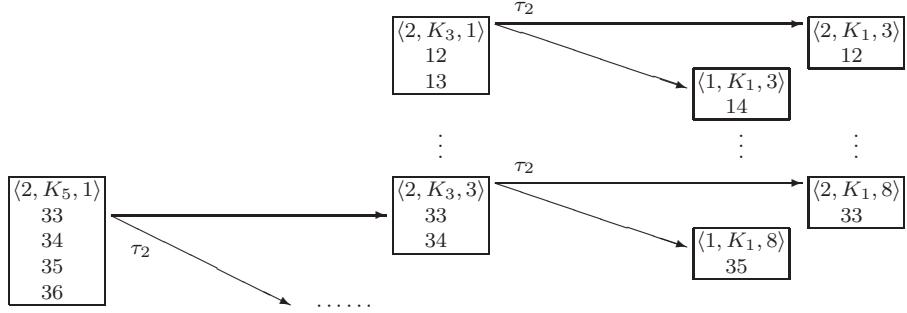
\centerline{Fig.3: An example of the graph embedding in the recursive structure $\mathcal{S}$.}

\vspace{0.2cm}

From Fig.3 we can see that: in the tree with root $\langle2,K_5,1\rangle$,
the branch from node $\langle2,K_3,3\rangle$ is the graph embedding of the tree with root $\langle2,K_3,1\rangle$.

\section{The number of repeated squares in $\mathbb{F}[1,n]$}

Denote $b(n):=\sharp\{(\omega,p):\omega_p\omega_{p+1}\triangleright\mathbb{F}[1,n]\}$
the number of squares ending at position $n$.
By the definition of $\langle i,K_m,p\rangle$, $b(n)$ is equal to the number of integer $n$ occurs in the recursive structure $\mathcal{S}$.
Thus we can calculate $b(n)$ by the property below.

\begin{property}[]\label{b}\ $b([1,2,3])=[0,0,0]$, $b(\Gamma_{2,1,1})=b([4])=[1]$, $b(\Gamma_{1,1,1})=b([5,6])=[0,1]$, $b(\Gamma_{2,2,1})=b([7,8])=[1,1]$, $b(\Gamma_{1,2,1})=b([9,10,11])=[1,1,2]$, for $m\geq3$,
\begin{equation*}
\begin{cases}
~b(\Gamma_{1,m,1})
=[b(\Gamma_{2,m-1,1}),b(\Gamma_{1,m-1,1})]+[0,\underbrace{1,\cdots,1}_{f_m-1}];\\
~b(\Gamma_{2,m,1})
=[b(\Gamma_{2,m-2,1}),b(\Gamma_{1,m-2,1})]
+[\underbrace{1,\cdots,1}_{f_{m-3}+1},\underbrace{0,\cdots,0}_{f_{m-2}-1}].
\end{cases}
\end{equation*}
\end{property}

The first few values of $b(n)$ are $b([1,2,3])=[0,0,0]$,
$b([4])=[1]$, $b([5,6])=[0,1]$,

$b([7,8])=[1,1]$, $b([9,10,11])=[1,1,2]$,
$b([12,13,14])=[2,1,1]$, $b([15,\cdots,19])=[1,2,2,2,3]$,

$b([20,\cdots,24])=[2,2,2,1,2]$, $b([25,\cdots,32])=[2,2,2,2,3,3,3,4]$.

\vspace{0.2cm}

For $m\geq3$, the immediately corollaries are
\begin{equation*}
\begin{cases}
\sum b(\Gamma_{1,m,1})=\sum b(\Gamma_{2,m-1,1})+\sum b(\Gamma_{1,m-1,1})+f_m-1;\\
\sum b(\Gamma_{2,m,1})=\sum b(\Gamma_{2,m-2,1})+\sum b(\Gamma_{1,m-2,1})+f_{m-3}+1.
\end{cases}
\end{equation*}


\begin{property}[]\label{P6.1} For $m\geq1$,
(1) $\sum b(\Gamma_{1,m,1})=\frac{2m+5}{5}f_{m}+\frac{2m-6}{5}f_{m-2}-1$,

(2) $\sum b(\Gamma_{2,m,1})=\frac{2m-2}{5}f_{m-1}+\frac{2m-3}{5}f_{m-3}+1$.
\end{property}

Since $\Gamma_{1,m,1}=[2f_{m+1}-1,\cdots,f_{m+3}-2]$ and
$\Gamma_{2,m,1}=[f_{m+2}-1,\cdots,2f_{m+1}-2]$, we have
$B(f_{m+3}-2)
=B(f_{m+2}-2)+\sum b(\Gamma_{2,m,1})+\sum b(\Gamma_{1,m,1})$ and
$B(2f_{m+1}-2)=B(f_{m+2}-2)+\sum b(\Gamma_{2,m,1})$. Thus by induction and Property \ref{P6.1},

\begin{property}[]\label{B}\
(1) $B(f_{m+3}-2)=\frac{2m-4}{5}f_{m+3}+\frac{2m}{5}f_{m+1}+4$ for $m\geq-1$.

(2) $B(2f_{m+1}-2)=\frac{4m-11}{5}f_{m+1}+\frac{4m-3}{5}f_{m-1}+5$ for $m\geq0$.
\end{property}





Property \ref{C6.4} can be proved by induction and Property \ref{b}.

\begin{property}\label{C6.4} $b(f_m-1)=\lfloor\frac{m-1}{2}\rfloor$, $b(f_m)=\lfloor\frac{m}{2}-1\rfloor$, $b(f_m-1)+b(f_m)=m-2$ for $m\geq2$.
\end{property}

\begin{remark}
Since $B(f_m)=B(f_{m}-2)+b(f_m-1)+b(f_m)$, $m\geq2$.
By Property \ref{B} and \ref{C6.4}, $$B(f_m)=\tfrac{4}{5}(m+1)f_{m}-\tfrac{2}{5}(m+7)f_{m-1}-4f_{m-2}+m+2
=\tfrac{4m-16}{5}f_{m}-\tfrac{2m-6}{5}f_{m-1}+m+2.$$
This is a known result of A.S.Fraenkel and J.Simpson\cite{FS2014}.
\end{remark}

Obversely we can calculate $B(n)$ by $B(n)=\sum_{i=4}^n b(i)$. But when $n$ is large, this method is complicated. Now we turn to give a fast algorithm.
For any $n\geq4$, let $m$ such that $f_m\leq n+1< f_{m+1}$. Since we already determine the expression of $B(f_m-2)$ and $B(2f_{m-1}-2)$ for $m\geq2$, in order to give a fast algorithm
of $B(n)$, we only need to calculate $\sum_{i=f_{m}-1}^nb(i)$ or $\sum_{i=2f_{m-1}-1}^nb(i)$. One method is calculating $b(n)$ by Property \ref{b}, the other method is using the corollaries as below.

\begin{corollary}\label{c1} For $n\geq4$, let $m$ such that $f_m\leq n+1\leq2f_{m-1}-1$, then $m\geq3$ and
\begin{equation*}
\sum_{i=f_{m}-1}^nb(i)
=\begin{cases}
\sum\limits_{i=f_{m-2}-1}^{n-f_{m-1}}b(i)+n-f_{m}+2,&n+1\leq f_{m}+f_{m-5}-1;\\
\sum\limits_{i=f_{m-2}+f_{m-5}-1}^{n-f_{m-1}}b(i)
+\frac{2m-5}{5}f_{m-5}+\frac{2m-11}{5}f_{m-7}+2,&otherwise.
\end{cases}
\end{equation*}
\end{corollary}

\begin{proof} By Property \ref{b},
when $f_m\leq n+1\leq f_{m}+f_{m-5}-1$,
$$\begin{array}{rl}
\sum\limits_{i=f_{m}-1}^nb(i)=\sum\limits_{i=f_{m-2}-1}^{n-f_{m-1}}[b(i)+1]
=\sum\limits_{i=f_{m-2}-1}^{n-f_{m-1}}b(i)+n-f_{m}+2.
\end{array}$$

When $f_{m}+f_{m-5}\leq n+1\leq2f_{m-1}-1$, $\sum\limits_{i=f_{m}-1}^nb(i)=
\sum\limits_{i=f_{m}-1}^{f_{m}+f_{m-5}-2}b(i)+\sum\limits_{i=f_{m}+f_{m-5}-1}^{n}b(i)$,
where
$$\left\{\begin{array}{rl}
\sum\limits_{i=f_{m}-1}^{f_{m}+f_{m-5}-2}b(i)
=&\sum\limits_{i=f_{m-2}-1}^{f_{m-2}+f_{m-5}-2}[b(i)+1]
=\sum b(\Gamma_{2,m-4,1})+f_{m-5}\\
=&\frac{2m-5}{5}f_{m-5}+\frac{2m-11}{5}f_{m-7}+1;\\
\sum\limits_{i=f_{m}+f_{m-5}-1}^{n}b(i)
=&\sum\limits_{i=f_{m-2}+f_{m-5}-1}^{n-f_{m-1}}b(i)+1.
\end{array}\right.$$
Thus $\sum\limits_{i=f_{m}-1}^nb(i)=
\sum\limits_{i=f_{m-2}+f_{m-5}-1}^{n-f_{m-1}}b(i)
+\frac{2m-5}{5}f_{m-5}+\frac{2m-11}{5}f_{m-7}+2$.
The conclusion holds.
\end{proof}

\begin{corollary}[]\label{c2}\
For $n\geq9$, let $m$ such that $2f_{m-1}\leq n+1\leq f_{m+1}-1$, then $m\geq4$ and
\begin{equation*}
\sum_{i=2f_{m-1}-1}^nb(i)
=\begin{cases}
\sum\limits_{i=f_{m-1}-1}^{n-f_{m-1}}b(i)+n-2f_{m-1}+1,
~~~~~~~~~~~~~~~~~~n+1\leq f_{m}+f_{m-2}-1;\\
\sum\limits_{i=2f_{m-2}-1}^{n-f_{m-1}}b(i)+n-2f_{m-1}
+\frac{2m-8}{5}f_{m-4}+\frac{2m-9}{5}f_{m-6}+2,~otherwise.
\end{cases}
\end{equation*}
\end{corollary}

\begin{proof} By Property \ref{b},
when $2f_{m-1}\leq n+1\leq 2f_{m-1}+f_{m-4}-1=f_{m}+f_{m-2}-1$,
$$\begin{array}{c}
\sum\limits_{i=2f_{m-1}-1}^nb(i)=\sum\limits_{i=f_{m-1}-1}^{n-f_{m-1}}b(i)+\sum\limits_{i=f_{m-1}}^{n-f_{m-1}}1
=\sum\limits_{i=f_{m-1}-1}^{n-f_{m-1}}b(i)+n-2f_{m-1}+1.
\end{array}$$
When $f_{m}+f_{m-2}\leq n+1\leq f_{m+1}-1$,
$\sum\limits_{i=2f_{m-1}-1}^nb(i)
=\sum\limits_{i=2f_{m-1}-1}^{f_{m}+f_{m-2}-2}b(i)
+\sum\limits_{i=f_{m}+f_{m-2}-1}^nb(i)$, where
$$\left\{\begin{array}{rl}
\sum\limits_{i=2f_{m-1}-1}^{f_{m}+f_{m-2}-2}b(i)
=&\sum\limits_{i=f_{m-1}-1}^{2f_{m-2}-2}b(i)+\sum\limits_{i=f_{m-1}}^{2f_{m-2}-2}1
=\sum\limits_{i=f_{m-1}-1}^{2f_{m-2}-2}b(i)+f_{m-4}-1\\
=&\sum b(\Gamma_{2,m-3,1})+f_{m-4}-1
=\frac{2m-3}{5}f_{m-4}+\frac{2m-9}{5}f_{m-6};\\
\sum\limits_{i=f_{m}+f_{m-2}-1}^nb(i)
=&\sum\limits_{i=2f_{m-2}-1}^{n-f_{m-1}}[b(i)+1]
=\sum\limits_{i=2f_{m-2}-1}^{n-f_{m-1}}b(i)+n-2f_{m-1}-f_{m-4}+2.
\end{array}\right.$$
Thus $\sum\limits_{i=2f_{m-1}-1}^nb(i)
=\sum\limits_{i=2f_{m-2}-1}^{n-f_{m-1}}b(i)+n-2f_{m-1}
+\frac{2m-8}{5}f_{m-4}+\frac{2m-9}{5}f_{m-6}+2$.
The conclusion holds.
\end{proof}

\noindent\textbf{Example.} One method to calculate $\sum_{i=20}^{23}b(i)$ is by Property \ref{b}.
Since $b(\Gamma_{2,4,1})=b([20,\cdots,24])=[2,2,2,1,2]$, $\sum_{i=20}^{23}b(i)=7$.
The other method is using Corollary \ref{c1} and \ref{c2}:
$$\begin{array}{c}
\sum\limits_{i=20}^{23}b(i)
=\sum\limits_{i=f_{4}+f_{1}-1}^{23-f_{5}}b(i)
+\frac{7}{5}f_{1}+\frac{1}{5}f_{-1}+2
=\sum\limits_{i=9}^{10}b(i)+5=7.
\end{array}$$

\begin{algorithm}[The number of repeated squares, $B(n)$]\

Step 1. For $n\leq3$, $B(n)=0$;
for $n\leq4$, find the $m$ such that $f_m\leq n+1< f_{m+1}$.

Step 2. Compare $n$ with $2f_{m-1}-1$.

(1) If $n<2f_{m-1}-1$, calculate $B(f_m-2)$ by Property \ref{B}; calculate $\sum_{i=f_{m}-1}^{n}b(i)$ by Property \ref{b} or by Corollary \ref{c1} and \ref{c2}.
Then $B(n)=B(f_{m}-2)+\sum_{i=f_{m}-1}^nb(i)$.

(2) If $n\geq2f_{m-1}-1$, calculate $B(2f_{m-1}-2)$ by Property \ref{B}; calculate $\sum_{i=2f_{m-1}-1}^nb(i)$ by Property \ref{b} or by Corollary \ref{c1} and \ref{c2}.
Then $B(n)=B(2f_{m-1}-2)+\sum_{i=2f_{m-1}-1}^nb(i)$.
\end{algorithm}

\begin{remark} When $m$ is large (resp. small), Corollary \ref{c1} and \ref{c2} (resp. Property \ref{b}) is faster.
\end{remark}

\noindent\textbf{Example.} We calculate $B(23)$. Since $f_{6}=21\leq 23+1<f_{7}=34$, $m=6$. Moreover $23<2f_{5}-1$.

By Property \ref{B}, $B(f_{6}-2)=B(19)=\frac{2}{5}f_{6}+\frac{6}{5}f_{4}+4=22$.
By Property \ref{b} or by Corollary \ref{c1} and \ref{c2}, $\sum_{i=20}^{23}b(i)=7$.
Thus $B(23)=B(19)+\sum_{i=20}^{23}b(i)=22+7=29$.

\section{Basic properties of cubes}

Let $\omega$ be a factor with kernel $K_m$, by an analogous argument as Section 3 and by Proposition 4.8 in \cite{HW2015-1},
$\omega_{p}\omega_{p+1}\omega_{p+2}\prec\mathbb{F}$ has only one case: $r_p(K_m)=r_{p+1}(K_m)=r_1(K_m)=K_mK_{m+1}$.
In this case, $|\omega|=f_{m+2}$. Moreover $2\leq i\leq f_{m+1}$ and $m\geq0$,
\begin{equation*}
\begin{split}
\omega\omega\omega=
&K_{m+1}[i,f_{m+1}] K_{m} K_{m+1}K_{m} K_{m+1} K_{m} K_{m+1}[1,i-1]\\
=&K_{m+2}[i,f_{m+2}] \underline{K_{m+3}} K_{m+2}[1,i+f_{m}-1]
=K_{m+6}[i+f_{m+3},i+f_{m+5}+f_{m}-1].
\end{split}
\end{equation*}
Since $K_{m+3}\prec\omega\omega\omega\prec K_{m+6}[2,f_{m+6}-1]$, by Property \ref{k2}, $Ker(\omega\omega\omega)=K_{m+3}$.

\begin{remark}
By the discussion above, we have:
all cubes in $\mathbb{F}$ are of length $3f_m$ for some $m\geq2$, and a cube of each such length occurs.
This is Theorem 8 in J.Shallit et al\cite{DMSS2014}.
\end{remark}

For $m\geq3$ and $p\geq1$, we define a set below:
$$\langle K_m,p\rangle:=
\{P(\omega\omega\omega,p):Ker(\omega\omega\omega)=K_m,|\omega|=f_{m-1},
\omega\omega\omega\prec\mathbb{F}\}.$$
Obviously it contains all cubes. By Property \ref{P} we have
\begin{equation*}
\begin{split}
\langle K_m,p\rangle=&\{P(\omega,p):
\omega=K_{m-1}[i,f_{m-1}] K_{m} K_{m-1}[1,i+f_{m-3}-1],2\leq i\leq f_{m-2}\}\\
=&\{P(K_m,p)+f_{m-3}+i-1,2\leq i\leq f_{m-2}\}\\
=&\{pf_{m+1}+\lfloor\phi p\rfloor f_{m}+2f_{m-1},\cdots,pf_{m+1}+\lfloor\phi p\rfloor f_{m}+f_{m+1}-2\}.
\end{split}
\end{equation*}


\begin{corollary}
$\sharp\langle K_m,p\rangle=f_{m-2}-1$ for $m\geq3$, $p\geq1$.
\end{corollary}

\section{The number of distinct cubes in $\mathbb{F}[1,n]$}

Denote $c(n):=\sharp\{\omega:\omega\omega\omega\triangleright\mathbb{F}[1,n],
\omega\omega\omega\not\!\prec\mathbb{F}[1,n-1]\}$. Obversely, $C(n)=\sum_{i=1}^n c(i)$.

\begin{property}[]
$\langle K_m,1\rangle=\{f_{m+1}+2f_{m-1},\cdots,2f_{m+1}-2\}$ for $m\geq3$.
\end{property}

Sets $\langle K_m,1\rangle$ are pairwise disjoint, and each set contains some consecutive integers. We get a chain
$\langle K_3,1\rangle=\{14\},\langle K_4,1\rangle=\{23,24\},\cdots,
\langle K_m,1\rangle,\cdots$.
So $c(n)=1$ iff $n\in\cup_{m\geq3}\langle K_m,1\rangle$.
Thus

\begin{property}[] $c(n)=1$ iff $n\in\cup_{m\geq3}\{f_{m+1}+2f_{m-1},\cdots,2f_{m+1}-2\}$.
\end{property}

By consider $C(2f_{m+1}-2)$ for $m\geq3$, we can give a fast algorithm of $C(n)$ for all $n\geq1$. Since $\sum_{i=-1}^mf_i= f_{m+2}-1$,
$C(2f_{m+1}-2)=\sum_{i=3}^{m}\sharp\langle K_i,1\rangle
=\sum_{i=3}^{m}(f_{i-2}-1)=f_{m}-m-1$.


\begin{theorem}[] For $n<14$, $C(n)=0$; for $n\geq14$, let $m$ s.t. $f_{m+1}+2f_{m-1}\leq n<f_{m+2}+2f_{m}-1$, then $m\geq3$ and
\begin{equation*}
C(n)=\begin{cases}
n-f_{m+1}-f_{m-1}-m+1,&n\leq 2f_{m+1}-2;\\
f_{m}-m-1,&otherwise.
\end{cases}
\end{equation*}
\end{theorem}

\begin{proof} When $2f_{m+1}-1\leq n\leq f_{m+2}+2f_{m}-1$, $c(n)=0$,
$C(n)=C(2f_{m+1}-2)=f_{m}-m-1$.

When $f_{m+1}+2f_{m-1}\leq n\leq 2f_{m+1}-2$, $c(n)=1$,
$$C(n)=C(f_{m+1}+2f_{m-1}-1)+n-f_{m+1}-2f_{m-1}+1.$$
Since $C(f_{m+1}+2f_{m-1}-1)=C(2f_{m}-2)=f_{m-1}-m$, $C(n)=n-f_{m+1}-f_{m-1}-m+1.$
Thus the conclusion holds.
\end{proof}

%
%

Since $2f_{m-2}-1\leq f_m\leq f_{m-1}+2f_{m-3}-1$ for $m\geq6$, we have

\begin{theorem}\
$C(f_m)=0$ for $m\leq5$, $C(f_{m})=f_{m-3}-m+2$ for $m\geq6$.
\end{theorem}

\section{The recursive structure of cubes}

In this section, we establish a recursive structure of cubes. Using it, we will count
the number of repeated cubes in $\mathbb{F}[1,n]$ (i.e. $D(n)$) in Section 10.

\begin{property}[]\label{R3}\ For $m\geq5$,
$\min\langle K_m,p\rangle-2=\max\langle K_{m-2},P(b,p)+1\rangle$.
\end{property}

\begin{proof} Since $P(b,p)=2p+\lfloor\phi p\rfloor$, $\lfloor\phi(2p+\lfloor\phi p\rfloor+1)\rfloor=p+\lfloor\phi p\rfloor$, for $m\geq3$, $f_{m-1}+f_{m-4}=2f_{m-2}$,
\begin{equation*}
\begin{split}
&\max\langle K_{m-2},P(b,p)+1\rangle+2=(2p+\lfloor\phi p\rfloor+1)f_{m-1}+\lfloor\phi (2p+\lfloor\phi p\rfloor+1)\rfloor f_{m-2}+f_{m-1}\\
=&(2p+\lfloor\phi p\rfloor+1)f_{m-1}+(p+\lfloor\phi p\rfloor)f_{m-2}+f_{m-1}
=pf_{m+1}+\lfloor\phi p\rfloor f_{m}+2f_{m-1}=\min\langle K_m,p\rangle.
\end{split}
\end{equation*}
This means $\max\langle K_{m-2},P(b,p)+1\rangle+2=\min\langle K_m,p\rangle$, so the
conclusion holds.
\end{proof}

By an analogous argument, we have

\begin{property}[]\label{R4}\ For $m\geq4$,
$\max\langle K_m,p\rangle+f_{m-4}+2=\min\langle K_{m-1},P(a,p)+1\rangle$.
\end{property}

In Property \ref{R3} and \ref{R4}, we establish the recursive relations for any
$\langle K_m,p\rangle$, $m\geq3$. Thus
we can define the recursive structure over $\{\langle K_m,p\rangle|~m\geq3,p\geq1\}$ denoted by $\mathcal{C}$. Each $\langle K_m,p\rangle$ is an element in $\mathcal{C}$.
The recursive structure $\mathcal{C}$ is a family of finite trees with roots $\langle K_m,1\rangle$ for all $m\geq3$; and with recursive relations:
\begin{equation*}
\begin{cases}
\tau_3\langle K_{m},p\rangle=\langle K_{m-2},P(b,p)+1\rangle\cup\langle K_{m-1},P(a,p)+1\rangle\text{ for }m\geq5;\\
\tau_4\langle K_4,p\rangle=\langle K_{m-1},P(a,p)+1\rangle.
\end{cases}
\end{equation*}
Since $\max\langle K_{m-2},P(b,p)+1\rangle<\min\langle K_{m-1},P(a,p)+1\rangle$, the ``$\cup$" is a disjoint union.

\begin{property}[]\
Each $\langle K_m,p\rangle$ belongs to the recursive structure $\mathcal{C}$, for $m\geq3$ and $p\geq1$.
\end{property}

\begin{proof} Each element $\langle K_m,1\rangle$ is root of a finite tree in $\mathcal{C}$.
For $p\geq1$,
\begin{equation*}
\begin{cases}
\langle K_m,P(a,p)+1\rangle\in\tau_3\langle K_{m+1},p\rangle~(m\geq4)
\text{ and }\langle K_3,P(a,p)+1\rangle\in\tau_4\langle K_{4},p\rangle;\\
\langle K_m,P(b,p)+1\rangle\in\tau_3\langle K_{m+2},p\rangle~(m\geq3).
\end{cases}
\end{equation*}
Since $\mathbb{N}=\{1\}\cup\{P(a,p)+1\}\cup\{P(b,p)+1\}$, the recursive structure $\mathcal{C}$ contains all $\langle K_m,p\rangle$.
\end{proof}

On the other hand, by the recursive relations $\tau_3$ and $\tau_4$, each element $\langle K_m,p\rangle$ has a unique position in $\mathcal{C}$.
Fig.4 show the finite tree in the recursive structure $\mathcal{C}$ with root $\langle K_6,1\rangle$.
\scriptsize
\setlength{\unitlength}{0.9mm}
\begin{center}
\begin{picture}(75,72)
\put(65,1){82}
\put(61.5,5){$\langle K_3,7\rangle$}
\put(61,0){\line(1,0){11}}
\put(61,0){\line(0,1){9}}
\put(72,0){\line(0,1){9}}
\put(61,9){\line(1,0){11}}
\put(65,29){69}
\put(61.5,33){$\langle K_3,6\rangle$}
\put(61,28){\line(1,0){11}}
\put(61,28){\line(0,1){9}}
\put(72,28){\line(0,1){9}}
\put(61,37){\line(1,0){11}}
\put(65,53){61}
\put(61.5,57){$\langle K_3,5\rangle$}
\put(61,52){\line(1,0){11}}
\put(61,52){\line(0,1){9}}
\put(72,52){\line(0,1){9}}
\put(61,61){\line(1,0){11}}
\put(45,5){79}
\put(45,9){78}
\put(41.5,13){$\langle K_4,4\rangle$}
\put(41,4){\line(1,0){11}}
\put(41,4){\line(0,1){13}}
\put(52,4){\line(0,1){13}}
\put(41,17){\line(1,0){11}}
\put(45,61){58}
\put(45,65){57}
\put(41.5,69){$\langle K_4,3\rangle$}
\put(41,60){\line(1,0){11}}
\put(41,60){\line(0,1){13}}
\put(52,60){\line(0,1){13}}
\put(41,73){\line(1,0){11}}
\put(25,13){74}
\put(25,17){73}
\put(25,21){72}
\put(25,25){71}
\put(21.5,29){$\langle K_5,2\rangle$}
\put(21,12){\line(1,0){11}}
\put(21,12){\line(0,1){21}}
\put(32,12){\line(0,1){21}}
\put(21,33){\line(1,0){11}}
\put(5,33){66}
\put(5,37){65}
\put(5,41){64}
\put(5,45){63}
\put(5,49){62}
\put(5,53){61}
\put(5,57){60}
\put(1.5,61){$\langle K_6,1\rangle$}
\put(1,32){\line(1,0){11}}
\put(1,32){\line(0,1){33}}
\put(12,32){\line(0,1){33}}
\put(1,65){\line(1,0){11}}
\put(14,45){$\tau_3$}
\put(13,40){\vector(1,1){27}}
\put(13,40){\vector(1,-2){7}}
\put(34,24){$\tau_3$}
\put(33,22){\vector(3,1){27}}
\put(33,22){\vector(1,-1){7}}
\put(54,10){$\tau_4$}
\put(53,10){\vector(1,-1){7}}
\put(54,65){$\tau_4$}
\put(53,65){\vector(1,-1){7}}
\end{picture}
\end{center}
\vspace{-0.2cm}
\normalsize
\centerline{Fig.4:  The finite tree in the recursive structure $\mathcal{C}$ with root $\langle K_6,1\rangle$.}

\begin{lemma}[] For $m\geq1$,
(1) $P(a,f_m-1)=f_{m+1}-2$, $\lfloor\phi (f_m-1)\rfloor=f_{m-1}-1$.

(2) $\lfloor\phi f_m\rfloor=f_{m-1}$ if $m$ is odd; $\lfloor\phi f_m\rfloor=f_{m-1}-1$
if $m$ is even.
(3) $P(b,f_{2m})=f_{2m+2}-1$.
\end{lemma}

\begin{proof} Denote by
$|\omega|_a$ (resp. $|\omega|_b$) the number of letter $a$ (resp. $b$) occurring in $\omega$.

(1) Since $|F_{m+1}|_a=f_{m}$, $aba\triangleright F_{2m}$ and $aab\triangleright F_{2m+1}$, we have $P(a,f_m-1)=f_{m+1}-2$.
On the other hand, by Corollary \ref{P1}, $P(a,f_m-1)=f_m-1+\lfloor\phi (f_m-1)\rfloor$.
Comparing the two expressions of $P(a,f_m-1)$, we have
$\lfloor\phi (f_m-1)\rfloor=f_{m-1}-1$ for $m\geq1$.

(2) By Corollary \ref{P1}, $P(a,f_m)=f_m+\lfloor\phi f_m\rfloor$. By the analogous argument in (1), we have: when $m$ is odd, $P(a,f_m)=f_{m+1}$, then $P(a,f_m)=f_{m+1}=f_m+\lfloor\phi f_m\rfloor\Rightarrow\lfloor\phi f_m\rfloor=f_{m-1}$;
when $m$ is even, $P(a,f_m)=f_{m+1}-1$, then $P(a,f_m)=f_{m+1}-1=f_m+\lfloor\phi f_m\rfloor\Rightarrow\lfloor\phi f_m\rfloor=f_{m-1}-1$.

(3) Since $|F_m|_b=f_{m-2}$, $aba\triangleright F_{2m}$, we have $P(b,f_{2m})=f_{2m+2}-1$ for $m\geq1$.
\end{proof}

\begin{lemma}[] $f_{m}f_{k}+f_{m-1}f_{k-1}=f_{m+k+1}$ for $m,k\geq-1$.
\end{lemma}

\begin{proof} Since $f_{m}f_{k}+f_{m-1}f_{k-1}=f_{m}(f_{k-1}+f_{k-2})+f_{m-1}f_{k-1}
=f_{m}f_{k-2}+(f_{m}+f_{m-1})f_{k-1}$,
using it repeatedly, $f_{m}f_{k}+f_{m-1}f_{k-1}=f_{m}f_{k-2}+f_{m+1}f_{k-1}=\cdots=f_{m+k-1}f_{-1}+f_{m+k}f_{0}=f_{m+k+1}$.
\end{proof}

For $m\geq3$, we define the vectors
$\Gamma_{m}:=[f_{m+2}-1,\cdots,f_{m+3}-2]$, then

\begin{property}[]\label{P9.6}\
The finite tree with root $\langle K_m,1\rangle$ belongs to $\Gamma_{m}$ for $m\geq3$.
\end{property}

\begin{proof} (1) Since $P(a,f_m-1)=f_{m+1}-2$, the maximal of the recursive structure from $\langle K_{m},1\rangle$ is
$$\max\{\max\langle K_{m},1\rangle,\max\langle K_{m-1},f_2-1\rangle,\max\langle K_{m-2},f_3-1\rangle,\cdots,\max\langle K_3,f_{m-2}-1\rangle\}$$
By Property \ref{R4}, $\max\langle K_{m},p\rangle<\min\langle K_{m-1},P(a,p)+1\rangle$,
so $\max\langle K_{m-i},f_{i+1}-1\rangle$ is strictly increasing for $0\leq i\leq m-3$.
Thus the maximal integer in the tree is $\max\langle K_3,f_{m-2}-1\rangle$.
\begin{equation*}
\begin{split}
&\max\langle K_3,f_{m-2}-1\rangle=(f_{m-2}-1)f_{4}+\lfloor\phi(f_{m-2}-1)\rfloor f_{3}+f_{4}-2\\
=&(f_{m-2}-1)f_4+(f_{m-3}-1)f_3+6=f_{m-2}f_4+f_{m-3}f_3-7=f_{m+3}-7<\max\Gamma_m.
\end{split}
\end{equation*}

(2) Similarly, since $P(b,f_{2m})=f_{2m+2}-1$,
$\min\langle K_{m-2i},f_{2i}\rangle$ is strictly decreasing for $0\leq i\leq[\frac{m-4}{2}]$. So the minimal integer in the tree is
\begin{equation*}
\min\{\min\langle K_{m},1\rangle,\min\langle K_{m-2},f_2\rangle,\min\langle K_{m-4},f_4\rangle,\cdots\}=
\begin{cases}
\min\langle K_4,f_{m-4}\rangle& \text{if $m$ is even;}\\
\min\langle K_3,f_{m-3}\rangle& \text{if $m$ is odd.}
\end{cases}
\end{equation*}

When $m$ is even, $\min\langle K_4,f_{m-4}\rangle=f_{m-4}f_{5}+\lfloor\phi f_{m-4}\rfloor f_{4}+2f_{3}$, $\lfloor\phi f_{m-4}\rfloor=f_{m-5}-1$, so
$$\min\langle K_4,f_{m-4}\rangle
=f_{m-4}f_{5}+(f_{m-5}-1)f_{4}+2f_{3}=f_{m+2}+2>\min\Gamma_m.$$

When $m$ is odd, $\min\langle K_3,f_{m-3}\rangle=f_{m-3}f_{4}+\lfloor\phi f_{m-3}\rfloor f_{3}+2f_{2}$, $\lfloor\phi f_{m-3}\rfloor=f_{m-4}-1$, so
$$\min\langle K_3,f_{m-3}\rangle
=f_{m-3}f_{4}+(f_{m-4}-1)f_{3}+2f_{2}
=f_{m+2}+1>\min\Gamma_m.$$

In each case, the minimal integer in the tree is larger than $\min\Gamma_m$, so the conclusion holds.
\end{proof}

By Property \ref{P9.6} and the definition of $\Gamma_m$, the finite trees in recursive structure $\mathcal{C}$ with different roots $\langle K_{m},1\rangle$ are disjoint.

\section{The number of repeated cubes in $\mathbb{F}[1,n]$}

Denote $d(n):=\sharp\{(\omega,p):\omega_p\omega_{p+1}\omega_{p+2}\triangleright\mathbb{F}[1,n]\}$, the number of cubes ending at position $n$. Obversely, $D(n)=\sum_{i=1}^n d(i)$.
By the definition of $\langle K_m,p\rangle$, $d(n)$ is equal to the number of integer $n$ occurs in the recursive structure $\mathcal{C}$.
Thus we can calculate $d(n)$ by the property below.

\begin{property}[]\label{d} For $m\geq3$,
\begin{equation*}
\begin{split}
&d([f_{m+4}-1,\cdots,f_{m+5}-2])\\
=&d([f_{m+2}-1,\cdots,f_{m+3}-2,f_{m+3}-1,\cdots,f_{m+4}-2])
+[\underbrace{0,\cdots,0}_{f_{m-1}+1},\underbrace{1,\cdots,1}_{f_m-1},\underbrace{0,\cdots,0}_{f_{m+2}}].
\end{split}
\end{equation*}
\end{property}

The first few values of $d(n)$ are
$d([f_{5}-1,\cdots,f_{6}-2])=[d(12),\cdots,d(19)]=[0,0,1,0,0,0,0,0]$,

$d([f_{6}-1,\cdots,f_{7}-2])=[d(20),\cdots,d(32)]=[0,0,0,1,1,0,0,1,0,0,0,0,0]$,

$d([f_{7}-1,\cdots,f_{8}-2])=[d(33),\cdots,d(53)]=[0,0,1,0,1,1,1,1,0,0,0,1,1,0,0,1,0,0,0,0,0]$.

\vspace{0.2cm}

By Property \ref{d}, $\sum d(\Gamma_{m+2})=\sum d(\Gamma_{m})+\sum d(\Gamma_{m+1})+f_{m-2}-1$. By induction, we have

\begin{lemma}[]
$\sum\limits_{f_{m+2}-1}^{f_{m+3}-2} d(n)
=\frac{m-5}{5}f_{m}+\frac{m+2}{5}f_{m-2}+1$ for $m\geq3$.
\end{lemma}



By the definition of $D(n)$, $D(f_{m+1}-2)=D(f_{m}-2)+\sum_{f_{m}-1}^{f_{m+1}-2} d(n)$. By induction, we have

\begin{property}[]\label{D}
$D(f_m-2)
=\frac{m-11}{5}f_{m-1}+\frac{m+1}{5}f_{m-3}+m+1$ for $m\geq6$.
\end{property}


By Property \ref{d}, we get $d(f_m-1)=d(f_m)=0$ easily by induction, thus
$D(f_m)=D(f_m-2)$.

\begin{theorem}[]\label{T10.4}
$D(f_m)
=\frac{m-11}{5}f_{m-1}+\frac{m+1}{5}f_{m-3}+m+1$ for $m\geq6$.
\end{theorem}

\begin{remark}
Theorem 59 in \cite{DMSS2014} shows the number of cube occurrences in $F_m$ as
$$D(f_m)=[d_1(m+2)+d_2]\alpha^{m+2}+[d_3(m+2)+d_4]\beta^{m+2}+m+1.$$
where $\alpha=\frac{1+\sqrt{5}}{2}$, $\beta=\frac{1-\sqrt{5}}{2}$,
$d_1=\frac{3-\sqrt{5}}{10}$, $d_2=\frac{17}{50}\sqrt{5}-\frac{3}{2}$,
$d_3=\frac{3+\sqrt{5}}{10}$, $d_4=-\frac{17}{50}\sqrt{5}-\frac{3}{2}$.
Since
$$\begin{array}{c}
f_m=\frac{\alpha^{m+2}-\beta^{m+2}}{\alpha-\beta},~
\alpha\beta=-1,~\frac{1}{\alpha}=\frac{\sqrt{5}-1}{2},~
\frac{1}{\beta}=\frac{-1-\sqrt{5}}{2},~ (\frac{1}{\alpha})^3=\sqrt{5}-2,~(\frac{1}{\beta})^3=-\sqrt{5}-2,
\end{array}$$
we can prove the two expressions are same.
By our expression in Theorem \ref{T10.4},
$$\begin{array}{rl}
&D(f_m)-m-1=\frac{m-11}{5}\times\frac{\alpha^{m+1}-\beta^{m+1}}{\alpha-\beta}
+\frac{m+1}{5}\times\frac{\alpha^{m-1}-\beta^{m-1}}{\alpha-\beta}\\
=&\frac{m-11}{5\sqrt{5}}\times\left(\frac{\sqrt{5}-1}{2}\alpha^{m+2}
+\frac{1+\sqrt{5}}{2}\beta^{m+2}\right)
+\frac{m+1}{5\sqrt{5}}\times\left((\sqrt{5}-2)\alpha^{m+2}+(\sqrt{5}+2)\beta^{m+2}\right)\\
=&[\frac{m-11}{5\sqrt{5}}\times\frac{\sqrt{5}-1}{2}+\frac{m+1}{5\sqrt{5}}(\sqrt{5}-2)]\alpha^{m+2}
+[\frac{m-11}{5\sqrt{5}}\times\frac{1+\sqrt{5}}{2}+\frac{m+1}{5\sqrt{5}}(\sqrt{5}+2)]\beta^{m+2}\\
=&[\frac{3-\sqrt{5}}{10}m+\frac{7\sqrt{5}-45}{50}]\alpha^{m+2}
+[\frac{3+\sqrt{5}}{10}m+\frac{-7\sqrt{5}-45}{50}]\beta^{m+2}.
\end{array}$$

By J.shallit's expression in \cite{DMSS2014},
$$\begin{array}{c}
D(f_m)-m-1=[\frac{3-\sqrt{5}}{10}m+\frac{3-\sqrt{5}}{5}+\frac{17}{50}\sqrt{5}-\frac{3}{2}]\alpha^{m+2}
+[\frac{3+\sqrt{5}}{10}m+\frac{3+\sqrt{5}}{5}-\frac{17}{50}\sqrt{5}-\frac{3}{2}]\beta^{m+2}.
\end{array}$$
Comparing the coefficients of $m\alpha^{m+2}$, $\alpha^{m+2}$, $m\beta^{m+2}$ and $\beta^{m+2}$, we have the two expressions are same.
\end{remark}

For any $n\geq12$, let $m$ such that $f_m\leq n+1< f_{m+1}$. Since we already determine the expression of $D(f_m-2)$, in order to give a fast algorithm
of $D(n)$, we only need to calculate $\sum_{i=f_{m}-1}^nd(i)$. One method is calculating $d(n)$ by Property \ref{d}, the other method is using the corollaries as below.

\begin{corollary}[]\label{c3}\
For $n\geq12$, let $m$ such that $f_m\leq n+1< f_{m+1}$, then $m\geq5$ and
\begin{equation*}
\sum_{i=f_{m}-1}^nd(i)
=\begin{cases}
\sum\limits_{i=f_{m-2}-1}^{n-f_{m-1}}d(i),&f_{m}\leq n+1\leq f_{m}+f_{m-5};\\
\sum\limits_{i=f_{m-2}+1}^{n-f_{m-1}}d(i)+n-f_{m}-f_{m-5}+1,&f_{m}+f_{m-5}+1\leq n+1\leq f_{m}+f_{m-3}-1;\\
\sum\limits_{i=f_{m-1}-1}^{n-f_{m-1}}d(i)+\frac{m-4}{5}f_{m-4}+\frac{m-2}{5}f_{m-6},&f_{m}+f_{m-3}\leq n+1< f_{m+1}.
\end{cases}
\end{equation*}
\end{corollary}

\begin{proof} By Property \ref{d},
when $f_{m}\leq n+1\leq f_{m}+f_{m-5}$,
$\sum\limits_{i=f_{m}-1}^nd(i)
=\sum\limits_{i=f_{m-2}-1}^{n-f_{m-1}}d(i)$.

When $f_{m}+f_{m-5}+1\leq n+1\leq f_{m}+f_{m-3}-1$,
$$\begin{array}{rl}
&\sum\limits_{i=f_{m}-1}^nd(i)
=\sum\limits_{i=f_{m}-1}^{f_{m}+f_{m-5}-1}d(i)
+\sum\limits_{i=f_{m}+f_{m-5}}^{n}d(i)
=\sum\limits_{i=f_{m-2}+1}^{n-f_{m-1}}d(i)
+\sum\limits_{i=f_{m-2}+f_{m-5}}^{n-f_{m-1}}1\\
=&\sum\limits_{i=f_{m-2}+1}^{n-f_{m-1}}d(i)+n-f_{m}-f_{m-5}+1.
\end{array}$$

When $f_{m}+f_{m-3}\leq n+1< f_{m+1}$,
$$\begin{array}{rl}
&\sum\limits_{i=f_{m}-1}^nd(i)
=\sum\limits_{i=f_{m}-1}^{f_{m}+f_{m-5}-1}d(i)
+\sum\limits_{i=f_{m}+f_{m-5}}^{f_{m}+f_{m-3}-2}d(i)
+\sum\limits_{i=f_{m}+f_{m-3}-1}^{n}d(i)\\
=&\sum\limits_{i=f_{m-2}-1}^{f_{m-1}-2}d(i)+f_{m-4}-1
+\sum\limits_{i=f_{m-1}-1}^{n-f_{m-1}}d(i)
=\sum\limits_{i=f_{m-1}-1}^{n-f_{m-1}}d(i)+\frac{m-4}{5}f_{m-4}+\frac{m-2}{5}f_{m-6}.
\end{array}$$

So the conclusion holds.
\end{proof}

\noindent\textbf{Example.}
Since
$[d(33),\cdots,d(53)]=[0,0,1,0,1,1,1,1,0,0,0,1,1,0,0,1,0,0,0,0,0]$,
using Property \ref{d},
we have $\sum_{i=33}^{48}d(i)=8$.
The other method is using Corollary \ref{c3}.
Since $f_{7}+f_{4}=42\leq n+1=49< f_{8}=55$,
$\sum_{i=33}^{48}d(i)
=\sum_{i=f_{6}-1}^{48-f_{6}}d(i)+\frac{3}{5}f_{3}+\frac{5}{5}f_{1}
=\sum_{i=20}^{27}d(i)+5
=\sum_{i=12}^{14}d(i)+7=8.$

\begin{algorithm}[The number of repeated cubes occurrences, $D(n)$]\

Step 1. For $n\leq11$, $D(n)=0$; for $n\leq12$, find the $m$ such that $f_m\leq n+1< f_{m+1}$.

Step 2. Calculate $D(f_m-2)$ by Property \ref{D}.

Step 3. Calculate $\sum_{i=f_{m}-1}^{n}d(i)$ by Property \ref{d} or by Corollary \ref{c3}.

Step 4. $D(n)=D(f_{m}-2)+\sum_{i=f_{m}-1}^{n}d(i)$.
\end{algorithm}

\noindent\textbf{Example.} We calculate $D(48)$. Since $f_{7}=34\leq 48+1< f_{8}=55$, $m=7$.

By Property \ref{D}, $D(32)=D(f_7-2)=\frac{-4}{5}f_{6}+\frac{8}{5}f_{4}+7+1=4$.

By Property \ref{d} or by Corollary \ref{c3}, $\sum_{i=33}^{48}d(i)=8$.
Thus $D(48)=D(32)+\sum_{i=33}^{48}d(i)=12$.

\vspace{0.5cm}

\noindent\textbf{\Large{Acknowledgments}}

\vspace{0.4cm}

The research is supported by the Grant NSF No.11431007, No.11271223 and No.11371210.

\end{CJK*}
\end{document}